\magnification=\magstep1
\input amstex
\documentstyle{amsppt}
\catcode`\@=11 \loadmathfont{rsfs}
\def\mycal{\mathfont@\rsfs}
\csname rsfs \endcsname \catcode`\@=\active

\NoBlackBoxes

\vsize=6.5in \topmatter
\title $L^2$-rigidity in von Neumann algebras.
\endtitle

\author JESSE PETERSON
\endauthor

\rightheadtext{$L^2$-rigidity}

\affil University of California, Los Angeles\endaffil

\address Math.Dept., UCLA, LA, CA 90095-155505\endaddress
\email jpete\@math.ucla.edu \endemail

\thanks \endthanks

\abstract  We introduce the notion of $L^2$-rigidity for von Neumann algebras,
a generalization of property (T) which can be viewed
as an analogue for the vanishing of 1-cohomology into the left regular representation of a group.  
We show that $L^2$-rigidity passes to normalizers and 
is satisfied by nonamenable II$_1$ factors
which are non-prime, have property $\Gamma$, or are weakly rigid.  
As a consequence we obtain that if $M$ is a free product of diffuse von Neumann algebras,
or if $M = L\Gamma$ where $\Gamma$ is a finitely generated group with $b_1^{(2)}(\Gamma) > 0$, then
any nonamenable regular subfactor of $M$ is prime and does not have properties $\Gamma$ or (T).
In particular this gives a new approach for showing primeness of all nonamenable subfactors of
a free group factor thus recovering a well known recent result of N. Ozawa.
\endabstract
\endtopmatter

\document
\heading 0. Introduction. \endheading

In their pioneering work of the 80's Connes and Jones ([C3], [CJ]) introduced the notion of 
property (T) (or rigidity) for II$_1$ factors by requiring that any sequence 
of subunital, subtracial completely positive maps which converge pointwise in $\| \cdot \|_2$
to the identity also converge uniformly in $\| \cdot \|_2$ to the identity on $(N)_1$.
This type of rigidity phenomenon (and it's relative version later introduced by Popa [P4]) has since led to the solution of
many old problems in von Neumann algebras and orbit equivalence ergodic theory ([C2], [IPP], [P3], [P4], [P5]).
In [Pe] it was shown that property (T) is equivalent to a vanishing 1-cohomology type result
for closable derivations into arbitrary Hilbert bimodules.  This equivalence is achieved in part
by using Sauvageot's results ([S1], [S3], [CiS]) which state that there is a bijective correspondence between
densely defined real closable derivations into Hilbert bimodules and semigroups of unital, tracial completely
positive maps.

For an inclusion of finite von Neumann algebras $(N \subset M)$ one cannot obtain such a characterization of relative 
property (T) introduced in [P4] (even if $N$ itself has property (T)) as there is no guarantee that a closed derivation 
$\delta$ on $M$ is even 
densely defined on $N$ much less inner.  However we will always have that the associated semigroup
will converge uniformly in $\| \cdot \|_2$ to id on $(N)_1$ and thus we may interpret this fact
as saying that $\delta$ ``vanishes'' on $N$.

In this paper we will use the techniques above to investigate closable derivations into the 
coarse correspondence $L^2(N) \overline \otimes L^2(N)$.
We will say that an inclusion of finite von Neumann algebra $(B \subset N)$ is $L^2$-rigid if 
all derivations which arise in this way ``vanish'' in the above sense on $B$, 
(see 3.1 for the precise definition).
Derivations into the coarse correspondence
appear naturally in the context of Voiculescu's nonmicrostates approach to free entropy [V], and 
also play a central role in studying the first $L^2$-Betti number of a von Neumann algebra as
introduced by Connes and Shlyakhtenko [CSh] (see also [T]).
This should be compared to the situation for groups where Bekka and Valette [BV] have shown that
for a finitely generated nonamenable group the first $L^2$-Betti number vanishes if
and only if the first cohomology group into the left regular representation vanishes.

We will show that given a nonamenable subfactor $Q \subset N$ and a densely defined real 
closable derivation into $(L^2(N) \overline \otimes L^2(N))^{\oplus \infty}$ then the derivation
must ``vanish'' on $Q' \cap N$.  Furthermore we will show that from the mixingness of
the coarse correspondence that if $Q' \cap N$ is diffuse then we further have that the 
derivation ``vanishes'' on $W^*(\Cal N_N(Q' \cap N))$.  Using a slight modification
of the above arguments using $N^\omega$ we will also show that if $N$ is a nonamenable 
factor which has property $\Gamma$ of Murray and von Neumann [MvN] then any derivation 
as above must ``vanish'' on $N$.  The main result is the following:

\proclaim{0.1. Theorem}  Let $N$ be a II$_1$ factor which is non-prime or has property
$\Gamma$, then $N$ is $L^2$-rigid.
\endproclaim

The above Theorem shows that $L^2$-rigidity is a very weak rigidity type phenomenon
(for instance $R \overline \otimes L\Bbb F_2$ is $L^2$-rigid even though it 
has Haagerup's compact approximation property [H]).  On the other hand we will
see that if $N$ is a free product of diffuse finite von Neumann algebras or 
if $N = L\Gamma$ where $\Gamma$ is a finitely generated group with $b_1^{(2)}(\Gamma) > 0$,
then $N$ is not $L^2$-rigid.
In [P1] Popa showed that for the uncountable free groups, their group factors are prime.
Using techniques from Voiculescu's free probability this was shown by Ge to also be the case for countable
free groups [Ge].  This was generalized to all free products of diffuse finite von Neumann algebras which
embed into $R^\omega$ by Jung [J].  From the above remarks we have the following:

\proclaim{0.2. Theorem}  Let $M$ be a free product of diffuse finite
von Neumann algebras or $M = L\Gamma$
where $\Gamma$ is a finitely generated group with $b_1^{(2)}(\Gamma) > 0$, then
any regular nonamenable subfactor of $M$ is prime and does not have properties $\Gamma$ or (T).
\endproclaim

Using techniques from $C^*$-algebra theory Ozawa
was able to show not just that the free group factors are prime but 
that in fact any subfactor of a free group factor is prime unless it is amenable [O1].
As an application of Theorem 0.1 we 
obtain a new approach to Ozawa's result using the fact that the free groups have
the ``$L^2$-Haagerup property'', i.e. there exist
proper cocycles into direct sums of the left regular representation.

\proclaim{0.3. Theorem}  Let $\Gamma$ be a countable discrete group such that there exists
a proper cocycle $b: \Gamma \rightarrow (\ell^2\Gamma)^{\oplus \infty}$,
(for example $\Gamma = \Bbb F_n$, $2 \leq n \leq \infty$).
Then any nonamenable subfactor of $L\Gamma$ is prime.
\endproclaim

It should be noted that although the above result gives a new proof of Ozawa's Theorem for the case
of the free groups, it is a quite different approach than in [O1].  Indeed, we use the fact that
$\Gamma$ has Haagerup's property in a crucial way.  Whereas in [O1] the above is shown for 
all hyperbolic groups, many of which have property (T).

In section 4 we investigate derivations which naturally appear in free products
of von Neumann algebras.  These derivations give rise to deformations by 
free products of multiples of the identity, thus we may extend the 
Kurosh type theorem in
[IPP] (Theorem 0.1) to include many von Neumann subalgebras which do not have relative property (T).
The first Kurosh type theorem in von Neumann algebras was obtained by Ozawa [O2] using 
$C^*$-algebra theory.

\proclaim{0.4. Theorem}  Let $M_1$ and $M_2$ be finite factors and let $M = M_1 * M_2$.
If $Q \subset M$ is a subfactor such that $Q' \cap M$ is a nonamenable factor,
or if $Q \subset M$ is a nonamenable subfactor with property $\Gamma$ and $Q' \cap M$ is a factor, 
then there exists $i \in \{ 1,2 \}$ and a unitary operator $u \in \Cal U(M)$ such that $uQu^* \subset M_i$.
\endproclaim

In section 5 we consider the case of a tensor product of II$_1$ factors 
$M = M_1 \overline \otimes \cdots \overline \otimes M_n$, such that each $M_i$ 
has a derivation into it's coarse correspondence which does not ``vanish''.  We show that
if $Q$ is a regular nonamenable subfactor then there exists a corner of $Q' \cap M$
which embeds into $M_i'$ for some $i \leq n$, where $M_i'$ is the von Neumann
subalgebra obtained by replacing $M_i$ with $\Bbb C$ in the above tensor product.
Ozawa and Popa [OP] gave examples of tensor products of von Neumann algebras
which have unique prime factorization.  Using the conjugacy results in [OP] 
we are able to give new examples of this type.

\proclaim{0.5. Theorem}  Let $M_i$ be nonamenable II$_1$ factors $1 \leq i \leq m$,
such that each $M_i$ is a non-trivial free product or $L\Gamma$ for some finitely generated group
$\Gamma$ with $b_1^{(2)}(\Gamma) > 0$, assume
$N_1 \overline \otimes \cdots \overline \otimes N_n = M_1 \overline \otimes \cdots \overline \otimes M_m$,
for some prime II$_1$ factors $N_1, \ldots, N_n$, then $n = m$ and there exist
$t_1, t_2, \ldots, t_m > 0$ with $t_1t_2 \cdots t_n = 1$ such that after permutation of indices 
and unitary conjugacy we have $N_k^{t_k} = M_k$, $\forall k \leq m$.
\endproclaim

This work was done while the author was at University of California, Los Angeles.  The author is very grateful
to the kind hospitality and stimulating environment at UCLA.  In particular the author would like to thank 
Adrian Ioana and Professor Sorin Popa for the many stimulating conversations regarding this work and 
for the insight the author has gained through the collaboration [IPP].  Also the author would like
to thank Kenley Jung for the useful discussions.

\heading 1. Preliminaries and notation. \endheading

Suppose $N$ is a finite von Neumann algebra with normal faithful trace $\tau$, 
$D(\delta) \subset N$ is a weakly dense $*$-subalgebra, 
$\Cal H$ is an $N$-$N$ Hilbert bimodule, and 
$\delta: D(\delta) \rightarrow \Cal H$ is a derivation
($\delta(xy) = x\delta(y) + \delta(x)y$, $\forall x,y \in D(\delta)$), which is closable
(as an unbounded operator from $L^2(N, \tau)$ to $\Cal H$), and real
($\langle \delta(x), y\delta(z) \rangle = \langle \delta(z^*) y^*, \delta(x^*) \rangle$, $\forall x,y,z \in D(\delta)$).  

It follows from [S1] and [DL] that $D(\overline \delta) \cap N$ is a $*$-subalgebra of $N$ and 
$\overline \delta|_{D(\overline \delta) \cap N}$ is again
a derivation. Let $\varDelta = \delta^* \overline \delta$, then $\varDelta$ is 
the generator of a completely Dirichlet form [S1].
Associated to $\varDelta$ are two natural deformations of $N$, the first is 
the completely positive semigroup (completely Markovian semigroup) 
$\{ \phi_t \}_{t > 0}$, each $\phi_t = \exp{(-t\varDelta)}$ is a c.p. map which is unital ($\phi_t(1) = 1$), tracial 
($\tau \circ \phi_t = \tau$), and positive ($\tau(\phi_t(x)x^*) \geq 0$, $\forall x \in N$),
moreover the semigroup property is satisfied ($\phi_{t+s} = \phi_t \circ \phi_s$, $\forall s,t>0$),
and $\forall x \in N$, $\| x - \phi_t(x) \|_2 \rightarrow 0$, as $t \rightarrow 0$.
The second deformation associated to $\varDelta$ is the deformation coming from resolvent maps
$\{ \eta_{\alpha} \}_{\alpha > 0}$, again each $\eta_{\alpha} = \alpha (\alpha + \varDelta)^{-1}$
is a unital, tracial, positive,
c.p. map such that $\forall x \in N$, $\| x - \eta_{\alpha}(x) \|_2 \rightarrow 0$, as $\alpha \rightarrow \infty$,
furthermore $\beta \eta_{\alpha} - \alpha \eta_{\beta} = (\beta - \alpha)\eta_{\alpha} \circ \eta_{\beta}$,
$\forall \alpha, \beta > 0$.

The relationship between these maps are as follows and can be found for example in [MR]:
$$
\varDelta = \lim_{t \rightarrow 0} {1 \over t} ({\text {id}} - \phi_t) 
        = \alpha( \eta_{\alpha}^{-1} - {\text {id}})
	= \lim_{\alpha \rightarrow \infty} \alpha({\text {id}} - \eta_{\alpha}),
$$
$$
\phi_t = {\text {exp}} (-t\varDelta)
	= \lim_{\alpha \rightarrow \infty} {\text {exp}}(-t\alpha({\text {id}} - \eta_{\alpha})),
$$
$$
\eta_{\alpha} = \alpha (\alpha + \varDelta)^{-1} 
	       = \alpha \int_0^{\infty} e^{-\alpha t}\phi_t dt.
$$

Note that we will use the same symbols $\varDelta, \phi_t$, and $\eta_\alpha$ for the maps on $N$ as well as
the corresponding extensions to $L^2(N, \tau)$.
Also note that $\eta_{\alpha}$ maps into the domain of $\varDelta$ 
and $\varDelta \circ \eta_{\alpha} = \alpha ({\text {id}} - \eta_{\alpha})$.  Furthermore
we have that Range$(\eta_\alpha) = D(\varDelta) \subset D(\overline \delta)$, 
$D(\varDelta^{1 \over 2}) = D(\overline \delta) = $Range$(\eta_\alpha^{1/2})$ and 
$\forall x \in D(\overline \delta)$, $\| \varDelta^{1 \over 2}(x) \|_2 = \| \delta(x) \|_2$.
	
If $B \subset N$ is a von Neumann subalgebra we will say that a deformation 
$\{ \Phi_\iota \}_\iota$ converges uniformly on $(B)_1$ if $\forall \varepsilon > 0$,
$\exists \iota_0$ such that $\forall \iota > \iota_0$, $b \in (B)_1$ we have
that $\| b - \Phi_\iota (b) \|_2 < \varepsilon$.  

\proclaim{1.1. Lemma}  Let $(N, \tau)$ be a finite von Neumann algebra, $B \subset N$ a von Neumann subalgebra,
and $\{ \phi_t \}_t$, $\{ \eta_\alpha \}_\alpha$ deformations as above.
The deformation $\{ \eta_\alpha \}_\alpha$ converges uniformly
on $(B)_1$ as $\alpha \rightarrow \infty$ if and only if the deformation $\{ \phi_t \}_t$ converges uniformly
on $(B)_1$ as $t \rightarrow 0$.
\endproclaim

\vskip .05in
\noindent
{\it Proof}.  Since $0 \leq \phi_t \leq id$, $\forall t > 0$ we have that
$\forall x \in N$, $t \mapsto \tau( (x - \phi_t(x)) x^*)$ is a non-negative valued function, also since 
$\tau ( (x - \phi_{t + s}(x)) x^*) = \tau( (x - \phi_t(x) ) x^* ) 
					+ \tau ( (\phi_{t/2} (x) - \phi_s( \phi_{t/2} (x))) \phi_{t/2}(x)^*) 
				\geq \tau( (x - \phi_t(x) ) x^* )$,
we have that $t \mapsto \tau( (x - \phi_t(x)) x^*)$ decreases to $0$ as 
$t \rightarrow 0$.  Hence if $\{ \phi_t \}_t$ does not converge uniformly on $(B)_1$ as $t \rightarrow 0$
then $\exists c_0 > 0$ such that $\forall t > 0$, $\exists x_t \in (B)_1$, such that
$\tau( (x_t - \phi_t(x_t)) x_t^* ) \geq c_0$.
Therefore $\tau( (x_t - \eta_{1/t}(x_t)) x_t^* ) = \int_0^\infty e^{-s} \tau( (x_t - \phi_{st}(x_t)) x_t^* ) ds
\geq \int_1^\infty e^{-s} c_0 ds \geq c_0(1 - e^{-1})$, thus $\{ \eta_\alpha \}_\alpha$ does not converge uniformly 
on $(B)_1$ as $\alpha \rightarrow \infty$.

Conversely if $\{ \phi_t \}_t$ does converge uniformly
on $(B)_1$ as $t \rightarrow 0$, then $\forall x \in (B)_1$ we have
$\| x - \eta_\alpha(x) \|_2 \leq \int_0^\infty e^s \| x - \phi_{s/\alpha} (x) \|_2 ds$
and since $\| x - \phi_t(x) \|_2 \leq 2$, $\forall x \in (B)_1$, $t > 0$ it follows that
$\{ \eta_\alpha \}_\alpha$ also converges uniformly 
on $(B)_1$ as $\alpha \rightarrow \infty$.
\hfill $\square$

\vskip .05in
Finally we mention that $\varDelta^{1 \over 2}$ also generates a completely Dirichlet form as is shown in [S3] by the formula:  
$\varDelta^{1 \over 2} = \pi^{-1} \int_0^{\infty} t^{-1/2}({\text {id}} - \eta_t) dt$.

\vskip .05in
\noindent
{\it Example 1}:  Suppose $\Gamma$ is a countable discrete group,
$\pi: \Gamma \rightarrow \Cal O(\Cal K)$ is an orthogonal representation, and
$b: \Gamma \rightarrow \Cal K$ is a 1-cocycle.  Then associated to this 
cocycle is a conditionally negative definite function $\psi$ given by
$\psi(\gamma) = \| b(\gamma) \|^2$, there is also a semigroup of positive
definite functions $\{ \varphi_t \}_t$ given by 
$\varphi_t(\gamma) = e^{-t\psi(\gamma)}$, and there is also the 
set of positive definite resolvents $\{ \chi_{\alpha} \}_{\alpha}$ given by
$\chi_{\alpha}(\gamma) = \alpha / (\alpha + \psi(\gamma))$.

Let $\Cal H = \Cal K \overline \otimes_{\Bbb R} L^2(L(\Gamma))$ and equip $\Cal H$ with the
$L(\Gamma)$ bimodule structure which satisfies 
$u_{\gamma} (\xi \otimes \xi') = \pi(\gamma)\xi \otimes u_{\gamma}\xi'$
and $(\xi \otimes \xi') u_{\gamma} = \xi \otimes \xi'u_{\gamma}$, 
$\forall \gamma \in \Gamma$, $\xi \in \Cal H$, $\xi' \in L^2(L\Gamma)$.
Let $\delta_b: \Bbb C \Gamma \rightarrow \Cal H$ be the derivation
which satisfies $\delta_b(u_{\gamma}) = b(\gamma) \otimes u_{\gamma}$,
$\forall \gamma \in \Gamma$,
then $\delta_b$ is a real closable derivation and so as described above
we can associated with $\delta_b$ the c.c.n. map $\varDelta$ along with 
the deformations $\{ \phi_t \}_t$ and $\{ \eta_{\alpha} \}_{\alpha}$.
It can be easily checked that we have the following relationships:
$$
\varDelta(u_{\gamma}) = \psi(\gamma) u_{\gamma}, \forall \gamma \in \Gamma,
$$
$$
\phi_t(u_{\gamma}) = \varphi_t(\gamma) u_{\gamma}, \forall \gamma \in \Gamma, t > 0,
$$
$$
\eta_{\alpha}(u_{\gamma}) = \chi_{\alpha}(\gamma) u_{\gamma}, \forall \gamma \in \Gamma, \alpha > 0.
$$
 
Note that in this case we have that if $\Lambda < \Gamma$ then the derivation ${\delta_b}_{|\Bbb C\Lambda}$ is inner
if and only if the cocycle $b_{|\Lambda}$ is inner if and only if the deformation $\{ \eta_\alpha \}_\alpha$
converges uniformly on $(L\Lambda)_1$. 
Note also that if $\Cal K$ is the left regular representation of $\Gamma$
then $\Cal H$ is the coarse correspondence for $L(\Gamma)$.

\vskip .05in
\noindent
{\it Example 2}:
Suppose $(M_1, \tau_1)$ and $(M_2, \tau_2)$ are finite diffuse von Neumann algebras, 
and let $(M, \tau) = (M_1 * M_2, \tau_1 * \tau_2)$.  If 
we let $\delta_i:M_1 *_{\text {Alg}} M_2 \rightarrow L^2(M) \otimes L^2(M)$ be the unique derivation
which satisfies $\delta_i(x) = x \otimes 1 - 1 \otimes x$, $\forall x \in M_i$ and 
$\delta_i(y) = 0$, $\forall y \in M_j$ where $j \not=i$.  Then it is easy to check that
$\delta_i$ defines a closable real derivation and a simple calculation 
(see for example Corollary 4.2 and the following remark in [Pe]) shows that the associated 
semigroups of c.p. maps are given by $\phi_s^1 = (e^{-2s}{\text {id}} + (1-e^{-2s})\tau) * {\text {id}}$,
and $\phi_s^2 = {\text {id}} * (e^{-2s}{\text {id}} + (1-e^{-2s})\tau)$.  
Hence we have that $\{ \phi_s^j \}_s$ does not converge uniformly on $(M)_1$ as $s \rightarrow 0$.

\heading 2. Approximation properties. \endheading

Throughout this section $\delta$ will be a real closable derivation, $\varDelta = \delta^* \overline \delta$
the corresponding generator of a completely Dirichlet form,
and also $\{ \eta_\alpha \}_\alpha$, and $\{ \phi_t \}_t$ will be the deformations described above.

\proclaim{2.1. Lemma} If $x, y, xy \in D(\varDelta)$, then 
$\| \varDelta(x)y + x \varDelta(y) - \varDelta(xy) \|_1 \leq 2 \| \delta(x) \| \| \delta(y) \|$.
\endproclaim

\vskip .05in
\noindent
{\it Proof}.  $\forall z \in D(\delta)$ such that $\| z \| \leq 1$ we have 
$$
| \tau( \varDelta(x)yz + x \varDelta(y)z - \varDelta(xy)z) | 
 = | \langle \delta(x), \delta(z^*y^*) \rangle + \langle \delta(y), \delta(x^*z^*) \rangle 
   - \langle \delta(xy), \delta(z^*) \rangle |
$$
$$
= | \langle \delta(x), \delta(z^*y^*) \rangle + \langle \delta(y), \delta(x^*z^*) \rangle 
  - \langle x\delta(y) + \delta(x)y, \delta(z^*)
\rangle |
$$
$$
= | \langle \delta(x), z^* \delta(y^*) \rangle + \langle \delta(y), \delta(x^*)z^* \rangle |
$$
$$
\leq \| \delta(x) \| \| z^* \delta(y^*) \| + \| \delta(y) \| \| \delta (x^*)z^* \|
\leq 2 \| \delta(x) \| \| \delta(y) \|.
$$
As $D(\delta)$ is weakly dense the result follows by applying Kaplansky's Theorem.
\hfill $\square$

\vskip .05in
\proclaim{2.2. Lemma} Let $\{ \eta_\alpha \}_\alpha$ be the deformation described above, $\forall \alpha > 0$, 
$\eta_{\alpha}^{1/2} = \pi^{-1} \int_0^{\infty} {t^{-1/2} \over 1 + t} \eta_{\alpha(1 + t) / t} dt$,
also $({\text {id}} - \eta_{\alpha})^{1/2} 
 = \pi^{-1} \int_0^{\infty} {t^{-1/2} \over 1 + t} ({\text {id}} - \eta_{t \alpha/(1+t)}) dt$.
\endproclaim

\vskip .05in
\noindent
{\it Proof}.  
$\forall \alpha > 0$, $t > 0$ we have:
$$
\eta_{\alpha} ( t + \eta_{\alpha})^{-1} 
  = \eta_{\alpha} ( (t(\alpha + \varDelta) + \alpha) (\alpha + \varDelta)^{-1} )^{-1}
$$
$$
  = {1 \over t} \eta_{\alpha} (\alpha + \varDelta) ( {\alpha(1+t) \over t} + \varDelta )^{-1}
  = {\alpha \over t} ( {\alpha(1+t) \over t} + \varDelta )^{-1}
  = {1 \over (1 + t)} \eta_{\alpha(1 +t)/t},
$$
Hence $\eta_{\alpha}^{1/2} = \pi^{-1} \int_0^{\infty} t^{-1/2} \eta_{\alpha} ( t + \eta_{\alpha})^{-1} dt
		= \pi^{-1} \int_0^{\infty} {t^{-1/2} \over 1 + t} \eta_{\alpha(1 + t) / t} dt$.

The formula for $({\text {id}} - \eta_{\alpha})^{1/2}$ is shown similarly.
\hfill $\square$

\vskip .05in
Since the range of $\eta_\alpha^{1/2}$ is the same as the domain of $\delta$ we may take the composition
$\delta \circ \eta_\alpha^{1/2}$ to obtain a bounded operator from $L^2(N, \tau)$ to $\Cal H$ whose norm is
no more than $(2\alpha)^{1/2}$.  
In fact $\alpha \| x - \eta_\alpha(x) \|_2^2 
		\leq \| \delta \circ \eta_\alpha^{1/2} (x) \|_2^2 
		= \alpha \tau ((x - \eta_\alpha(x))x^*)
		\leq \alpha \| x - \eta_\alpha(x) \|_2$, $\forall x \in N$.
It will be convenient therefore to use the following notation, we will let $\zeta_\alpha = \eta_\alpha^{1/2}$, and
we will let $\tilde \delta_\alpha = \alpha^{-1/2} \delta \circ \zeta_\alpha$.  
The next lemma shows that $\tilde \delta_\alpha$ is
almost a derivation.

\vskip .05in
\proclaim{2.3. Lemma} Using the same notation as above
$\forall F \subset (N)_1$, such that $\{ \eta_\alpha \}_\alpha$ converges uniformly on $F$ 
($F$ possibly infinite),
$\forall \varepsilon > 0$, $\exists \alpha_0 > 0$, 
such that $\forall \alpha \geq \alpha_0$ we have that 
$\| \tilde \delta_\alpha( ax ) - \zeta_\alpha(a) \tilde \delta_\alpha( x ) 
				- \tilde \delta_\alpha(a) \zeta_\alpha( x ) \|^2 < \varepsilon$,
$\forall a \in F$, $x \in (N)_1$.  
\endproclaim

\vskip .05in
\noindent
{\it Proof}.  We will prove the Lemma in two parts: $(a),(b)$.  First we show that the vectors 
$\tilde \delta_\alpha( ax )$ and 
$\zeta_\alpha(a) \tilde \delta_\alpha( x ) + \tilde \delta_\alpha(a) \zeta_\alpha( x )$
have approximately the same size, and then we show that the vectors have large 
inner product.  The main difficulty is that we may not apply the product rule to a vector
of the form $\alpha^{-1/2} \delta \circ \zeta_\alpha ( ax )$ and thus in order to estimate
the size on an inner product we must translate the expression in terms involving $\varDelta^{1 \over 2}$
and then use Lemmas 2.1 (with the generator $\varDelta^{1 \over 2}$) and 2.2 to estimate these expressions.  
However some care is involved 
here as Lemma 2.1 only gives an estimate in $\| \cdot \|_1$ and thus we must make sure that
when we apply 2.1 the term we are taking the inner product with is bounded in uniform norm.

Thus each of the parts above separate into three steps: $(1), (2), (3)$.  The first step we use the properties of
the derivation to set up the $\| \cdot \|_1$ estimate from Lemma 2.1, the second step we translate
to terms involving $\varDelta^{1 \over 2}$ and use 2.1, and then the third step we 
use Lemma 2.2 and then translate back into terms of the derivation to finish the estimate.

Let $F \subset (N)_1$ be given as above and let $\varepsilon > 0$.
It follows from Lemma 2.2 and Section 1.1.2 of [P4] that $\exists \alpha_0 > 0$ 
such that $\forall \alpha \geq \alpha_0$ we have 
$\| a - \eta_\alpha (a) \|_2 < (\varepsilon/64)^4$,
$\| a - \zeta_\alpha (a) \|_2 < \varepsilon/100$,
and $\| a ({\text {id}} - \eta_\alpha)^{1/2} (x) - ({\text {id}} - \eta_\alpha)^{1/2} (ax) \|_2
	\leq \pi^{-1} \int_0^{\infty} {t^{1/2} \over 1 + t} 
		\| a \eta_{t\alpha/(1 + t)} (x) - \eta_{t\alpha/(1 + t)} (ax) \|_2 < \varepsilon/100$,
$\forall a \in F$, $x \in (N)_1$.
Then by using the product rule for the derivation we have 
$$
| \alpha^{-1} \| \delta( \zeta_\alpha(a) \zeta_\alpha(x) ) \|^2 
	- \alpha^{-1} \langle \delta( \zeta_\alpha(x) ), 
		        \delta( \zeta_\alpha(a^*) \zeta_\alpha(a) \zeta_\alpha(x) ) \rangle |
\tag a.1
$$
$$
\leq 8 \| \tilde \delta_\alpha( a ) \| \leq 8 \| a - \eta_\alpha(a) \|_2^{1/2} 
< \varepsilon/8.
$$

By Lemma 2.1 we have
$$
| \alpha^{-1} \langle \varDelta^{1 \over 2} \circ \zeta_\alpha(x), 
		      \varDelta^{1 \over 2} ( \zeta_\alpha(a^*) \zeta_\alpha(a) \zeta_\alpha(x)  ) \rangle
\tag a.2
$$
$$
	- \alpha^{-1} \langle \varDelta^{1 \over 2} \circ \zeta_\alpha(x),
		      \zeta_\alpha(a^*) \zeta_\alpha(a) \varDelta^{1 \over 2} \circ \zeta_\alpha(x)  \rangle |
$$
$$
\leq 2 \alpha^{1/2} \| \varDelta^{1 \over 2} ( \zeta_\alpha(a^*) \zeta_\alpha(a) \zeta_\alpha(x) )
			- \zeta_\alpha(a^*) \zeta_\alpha(a) \varDelta^{1 \over 2} \circ \zeta_\alpha(x) \|_1
$$
$$
\leq 2 \alpha^{1/2} \| \varDelta^{1 \over 2} ( \zeta_\alpha(a^*) \zeta_\alpha(a) ) \|_1 
   + 4 \alpha^{1/4} \| \varDelta^{1 \over 4} ( \zeta_\alpha(a^*) \zeta_\alpha(a) ) \|_2
$$
$$
\leq 4 \| a - \eta_\alpha(a) \|_2^{1/2} + 8 \| a - \eta_\alpha(a) \|_2^{1/4}
< \varepsilon/4.
$$

Also from the assumptions above we have 
$$
\alpha^{-1} |  \| \zeta_\alpha(a) \varDelta^{1 \over 2} \circ \zeta_\alpha (x) \|_2^2 
		- \| \varDelta^{1 \over 2} \circ \zeta_\alpha (ax) \|_2^2  |
\tag a.3
$$
$$
\leq 4\alpha^{-1/2} \| \zeta_\alpha (a) \varDelta^{1 \over 2} \circ \zeta_\alpha (x) 
		     - \varDelta^{1 \over 2} \circ \zeta_\alpha (ax)\|_2
$$
$$
\leq 8 \| \zeta_\alpha (a) - a \|_2 
	+ 4 \| a ({\text {id}} - \eta_\alpha)^{1/2} (x) - ({\text {id}} - \eta_\alpha)^{1/2} (ax) \|_2
< \varepsilon/8.
$$

Hence by combining $(a.1), (a.2)$, and $(a.3)$ we have shown 
$$
|  \| \alpha^{-1/2} \delta( \zeta_\alpha(a) \zeta_\alpha(x) ) \|^2 -
	\| \tilde \delta_\alpha(ax) \|^2 | < \varepsilon/2.
\tag a
$$

Similarly by using the product rule we obtain that
$$
| \alpha^{-1} \langle \delta( \zeta_\alpha(a) \zeta_\alpha(x) ), \delta( \zeta_\alpha (ax) ) \rangle  
	- \alpha^{-1} \langle \delta( \zeta_\alpha(x) ), 
		        \delta( \zeta_\alpha(a^*) \zeta_\alpha(ax) ) \rangle |
\tag b.1
$$
$$
\leq 4 \| \tilde \delta_\alpha(a) \| \leq 4 \| a - \eta_\alpha(a) \|_2^{1/2}
<\varepsilon/ 16.
$$

Again by Lemma 2.1 we have
$$
|\alpha^{-1} \langle \varDelta^{1 \over 2} \circ \zeta_\alpha(x), 
		      \varDelta^{1 \over 2} ( \zeta_\alpha(a^*) \zeta_\alpha(ax) ) \rangle
\tag b.2
$$
$$
	- \alpha^{-1} \langle \varDelta^{1 \over 2} \circ \zeta_\alpha(x),
		      \zeta_\alpha(a^*) \varDelta^{1 \over 2} \circ \zeta_\alpha(ax) \rangle |
$$
$$
\leq 2 \alpha^{-1/2} \| \varDelta^{1 \over 2} ( \zeta_\alpha(a^*) \zeta_\alpha(ax) ) 
				- \zeta_\alpha(a^*) \varDelta^{1 \over 2} \circ \zeta_\alpha(ax) \|_1
$$
$$
\leq 2 \alpha^{-1/2} \| \varDelta^{1 \over 2} \circ \zeta_\alpha(a^*) \zeta_\alpha(ax) \|_1
   + 4 \alpha^{1/4} \| \varDelta^{1 \over 4} \circ \zeta_\alpha(a^*) \|_2
$$
$$
\leq 2 \| a - \eta_\alpha(a) \|_2^{1/2} + 4 \| a - \eta_\alpha(a) \|_2^{1/4}
< \varepsilon/ 8.
$$

Also from the assumptions above we have
$$
|\alpha^{-1} \langle  \zeta_\alpha(a) \varDelta^{1 \over 2} \circ \zeta_\alpha(x),
		      \varDelta^{1 \over 2} \circ \zeta_\alpha(ax) \rangle
	- \alpha^{-1} \| \varDelta^{1 \over 2} \circ \zeta_\alpha(ax) \|_2^2 |
\tag b.3
$$
$$
\leq 4 \| \zeta_\alpha (a) - a \|_2 
	+ 2 \| a({\text {id}} - \eta_\alpha)^{1/2} (x) - ({\text {id}} - \eta_\alpha)^{1/2} (ax) \|_2
<\varepsilon/ 16.
$$

Thus using $(b.1), (b.2)$, and $(b.3)$ we have 
$$
| \langle \alpha^{-1/2}\delta( \zeta_\alpha(a) \zeta_\alpha(x) ), \tilde \delta_\alpha(ax) \rangle
	- \| \tilde \delta_\alpha(ax) \|_2^2 | < \varepsilon/4.
\tag b
$$

Hence by $(a)$ and $(b)$ we have that 
$$
\| \tilde \delta_\alpha( ax ) - \zeta_\alpha(a) \tilde \delta_\alpha( x ) 
				- \tilde \delta_\alpha(a) \zeta_\alpha( x ) \|^2
$$
$$
= \| \tilde \delta_\alpha(ax) - \alpha^{-1/2}\delta( \zeta_\alpha(a) \zeta_\alpha(x) ) \|^2
$$
$$
= \| \tilde \delta_\alpha(ax) \|^2 
   - 2\Re{ \langle \alpha^{-1/2}\delta( \zeta_\alpha(a) \zeta_\alpha(x) ), \tilde \delta_\alpha (ax) \rangle }
   + \| \alpha^{-1/2} \delta( \zeta_\alpha(a) \zeta_\alpha(x) ) \|^2
$$
$$
\leq | \| \tilde \delta_\alpha(ax) \|^2 
	- \| \alpha^{-1/2} \delta( \zeta_\alpha(a) \zeta_\alpha(x) ) \|^2 |
$$
$$
   + 2| \langle \alpha^{-1/2} \delta( \zeta_\alpha(a) \zeta_\alpha(x) ), \tilde \delta_\alpha (ax) \rangle
	- \| \tilde \delta_\alpha(ax) \|^2 | < \varepsilon.
$$
\hfill $\square$

The vectors $\tilde \delta_\alpha (x)$ may not be right and left bounded and thus vectors of the form
$\zeta_\alpha(a) \tilde \delta_\alpha (x)$ and $a \tilde \delta_\alpha (x)$ may be far apart 
even if $a$ and $\zeta_\alpha (a)$ are close.  The next lemma will allow us to handle this type of 
situation by showing that the normal state associated with $\tilde \delta_\alpha (x)$
has nice properties once we compose it with $\zeta_\alpha$.

\proclaim{2.4. Lemma} Given $x \in N$, such that $\delta(x) \not= 0$, and
$\alpha > 0$ let $\psi_\alpha^x$ be the normal state given by 
$\psi_\alpha^x(y) = \| \tilde \delta_\alpha (x) \|^{-2} 
			\langle y \tilde \delta_\alpha (x), \tilde \delta_\alpha (x) \rangle$.

\vskip .05in
\noindent
(i).  For $y \in D(\delta)$, we have 
$\psi_\alpha^x(y) = {1 \over 2}\| \delta(\zeta_\alpha(x)) \|^{-2}
	(\langle \varDelta^{1 \over 2}( y \zeta_\alpha(x)),  \varDelta^{1 \over 2} (\zeta_\alpha(x)) \rangle
	    + \langle \varDelta^{1 \over 2}( \zeta_\alpha(x^*) y ), \varDelta^{1 \over 2} (\zeta_\alpha (x^*)) \rangle
	    - \langle \varDelta^{1 \over 2}( y ), \varDelta^{1 \over 2} ( \zeta_\alpha(x) \zeta_\alpha(x^*) ) \rangle)$.

\vskip .05in
\noindent
(ii).  $\forall \varepsilon > 0$, $F \subset N$ finite, $\exists \alpha_0 > 0$, 
such that $\forall x \in N$, $\alpha > \alpha_0$ we have that 
$| \psi_\alpha^x( \zeta_\alpha(az) ) - \psi_\alpha^x( \zeta_\alpha(a) \zeta_\alpha(z) ) | 
	< \varepsilon / \| x - \eta_\alpha(x) \|_2^2$,
$\forall a \in F$, $z \in (N)_1$.

\vskip .05in
\noindent
(iii).  $\forall x, y \in N$, $\alpha > 0$, 
$| \psi_\alpha^x \circ \zeta_\alpha (y) | \leq 20 \| y \|_2 / \| x - \eta_\alpha(x) \|_2^2$.
\endproclaim

\vskip .05in
\noindent
{\it Proof}. 
{\it (i)}.  This follows by using the Leibniz rule for $\delta$ as in Lemma 2.1. 

\vskip .05in
\noindent
{\it (ii)}.  Let $\varepsilon > 0$, by Lemma 2.3 $\exists \alpha_1 > 0$ 
such that $\forall \alpha > \alpha_1$, $a \in F$, $z \in (N)_1$ we have
$\alpha^{-1/2} \| \varDelta^{1 \over 2} ( \zeta_\alpha(az) - \zeta_\alpha(a) \zeta_\alpha(z) ) \|_2 
			< {1 \over 9}(\varepsilon/16)^2$.
By Lemma 2.2 and 1.1.2 in [P4] let $\alpha_0 \geq \alpha_1$ such that $\forall \alpha \geq \alpha_0$ we have
$\| \zeta_\alpha(az) - a\zeta_\alpha(z) \|_1 < {1 \over 3} (\varepsilon/16)$, $\forall a \in F$, $z \in (N)_1$.

Then by Lemma 2.1 we have that $\forall a \in F$, $z \in (N)_1$, $\alpha \geq \alpha_0$,
$$
\alpha^{-1} | \langle \varDelta^{1 \over 2} ( ( \zeta_\alpha(az) 
					- \zeta_\alpha(a) \zeta_\alpha(z) ) \zeta_\alpha (x) ), 
		\varDelta^{1 \over 2} ( \zeta_\alpha (x) ) \rangle |
$$
$$
\leq 2\alpha^{-1/2} \| \varDelta^{1 \over 2} ( ( \zeta_\alpha(az) - \zeta_\alpha(a) \zeta_\alpha(z) ) 
					\zeta_\alpha (x) ) \|_1
$$
$$
\leq 4 \| \zeta_\alpha(az) - \zeta_\alpha(a) \zeta_\alpha(z) \|_1 
	+ 2\alpha^{-1/2} \| \varDelta^{1 \over 2} ( \zeta_\alpha(az) - \zeta_\alpha(a) \zeta_\alpha(z) ) \|_1
$$
$$
	+ 16\alpha^{-1/4} \| \varDelta^{1 \over 2} ( \zeta_\alpha(az) - \zeta_\alpha(a) \zeta_\alpha(z) ) \|_2^{1/2}
< 2\varepsilon/3.
$$

Similarly we have that
$$
\alpha^{-1} | \langle \varDelta^{1 \over 2} ( \zeta_\alpha (x^*)( \zeta_\alpha(az) 
								- \zeta_\alpha(a) \zeta_\alpha(z) ) ), 
		\varDelta^{1 \over 2} ( \zeta_\alpha (x^*) ) \rangle | < 2\varepsilon/3.
$$

Also
$$
\alpha^{-1} | \langle \varDelta^{1 \over 2} ( \zeta_\alpha(az) - \zeta_\alpha(a) \zeta_\alpha(z) ), 
		\varDelta^{1 \over 2} ( \zeta_\alpha (x) \zeta_\alpha (x^*) ) \rangle |
$$
$$
\leq 4 \alpha^{-1/2} \| \varDelta^{1 \over 2} ( \zeta_\alpha(az) - \zeta_\alpha(a) \zeta_\alpha(z) ) \|_2 < 2\varepsilon/3.
$$

Hence by (i) and the triangle inequality we have that
$| \psi_\alpha^x( \zeta_\alpha(az) ) - \psi_\alpha^x( \zeta_\alpha(a) \zeta_\alpha(z) ) | 
	< \varepsilon / \| \tilde \delta_\alpha( x ) \|^2 < \varepsilon / \| x - \eta_\alpha(x) \|_2^2$.

\vskip .05in
\noindent
{\it (iii)}.  Let $x,y \in N$, $\alpha > 0$, then

$$
\alpha^{-1} | \langle \varDelta^{1 \over 2} ( \zeta_\alpha(y) \zeta_\alpha (x) ), 
		\varDelta^{1 \over 2} ( \zeta_\alpha (x) ) \rangle |
$$
$$
\leq 2\alpha^{-1/2} \| \varDelta^{1 \over 2} ( \zeta_\alpha(y) \zeta_\alpha (x)) \|_1
$$
$$
\leq 4 \| \zeta_\alpha(y) \|_1 + 2\alpha^{-1/2} \| \varDelta^{1 \over 2} ( \zeta_\alpha(y) )\|_1 
			+ 4 \alpha^{-1/4} \| \varDelta^{1 \over 4}(\zeta_\alpha(y)) \|_2 \leq 16 \| y \|_2. 
$$

Similarly we have that
$$
\alpha^{-1} | \langle \varDelta^{1 \over 2} ( \zeta_\alpha (x^*) \zeta_\alpha(y) ), 
		\varDelta^{1 \over 2} ( \zeta_\alpha (x^*) ) \rangle | \leq 16 \| y \|_2.
$$

Also
$$
\alpha^{-1} | \langle \varDelta^{1 \over 2} ( \zeta_\alpha(y) ), 
		\varDelta^{1 \over 2} ( \zeta_\alpha (x) \zeta_\alpha (x^*) ) \rangle |
$$
$$
\leq 4 \alpha^{-1/2} \| \varDelta^{1 \over 2} (\zeta_\alpha(y) ) \|_2 \leq 8 \| y \|_2.
$$

Hence just as above we have that 
$| \psi_\alpha^x \circ \zeta_\alpha(y) | \leq 20 \| y \|_2 / \| x - \eta_\alpha(x) \|_2^2$.
\hfill $\square$

\heading 3. $L^2$-rigidity. \endheading

\proclaim{3.1. Definition} Let $N$ be a finite von Neumann algebra with trace $\tau$, 
if $M$ is a finite von Neumann algebra with trace $\tau'$ such that $N \subset M$, $\tau'|_N = \tau$,
and $\delta$ is a densely defined real closable derivation on $M$ into
$( L^2(M, \tau') \overline \otimes L^2(M, \tau') )^{\oplus \infty}$
then we say that the associated deformation 
$\{ \eta_\alpha \}_\alpha$ is an $L^2$-deformation for $N$.

If $B \subset N$ is a von Neumann subalgebra, the inclusion $(B \subset N)$ is $L^2$-rigid 
(or $B$ is an $L^2$-rigid subalgebra of $N$) if any $L^2$-deformation for $N$
converges uniformly on $(B)_1$.  
We will say that $N$ is $L^2$-rigid if the inclusion $(N \subset N)$ is $L^2$-rigid.
\endproclaim

\vskip .05in

\vskip .05in
\noindent
{\bf 3.2. Remarks}.  $1.$  It follows trivially that if $(B \subset N)$ is a rigid inclusion in the sense of [P4] then
$(B \subset N)$ is $L^2$-rigid.

\vskip .05in
\noindent
$2.$  By the definition it follows that if $M$ is a finite von Neumann algebra with 
normal faithful trace $\tau$
and $B \subset N \subset M$ are von Neumann subalgebras, then $(B \subset M)$ is $L^2$-rigid if 
$(B \subset N)$ is $L^2$-rigid.  

\vskip .05in
\noindent
$3.$  If $\Gamma$ is a discrete group such that $H^1(\Gamma, \ell^2\Gamma) \not= \{ 0 \}$
then from Example 1 in Section 1 we have that $L\Gamma$ is not $L^2$-rigid.
Also if $(M_1, \tau_1)$ and $(M_2, \tau_2)$ are finite diffuse von Neumann algebras 
then from Example 2 in Section 1 we have that
$(M_1 * M_2, \tau_1 * \tau_2)$ is not $L^2$-rigid.

\vskip .05in
\noindent
$4.$  If $\Gamma$ is a countable discrete group which has a proper cocycle
$b:\Gamma \rightarrow (\ell^2\Gamma)^{\oplus \infty}$ 
(for instance $\Gamma = \Bbb F_n, 1 \leq n \leq \infty$)
then $L\Gamma$ has no diffuse $L^2$-rigid von Neumann subalgebra.  Indeed if $\{ \eta_\alpha \}_\alpha$ 
is the associated deformation then $\eta_\alpha \in \Cal K(L^2(L\Gamma))$, $\forall \alpha > 0$ and 
thus if $B \subset L\Gamma$ is a von Neumann subalgebra
such that $\forall \varepsilon > 0$,
$\exists \alpha_0 > 0$ such that $\forall \alpha > \alpha_0$, $x \in B_1$ we have 
$\| x - \eta_{\alpha}(x) \|_2 < \varepsilon$ then we must have that $B$ cannot be diffuse 
(see for example Theorem 5.4 in [P4]).

\vskip .1in
Suppose $\Gamma = \Gamma_1 \times \Gamma_2$ where $\Gamma_1$ is infinite and $\Gamma_2$ is nonamenable, let us now 
sketch a simple proof that $H^1(\Gamma, \ell^2\Gamma) = \{ 0 \}$ (see also Corollary 10 in [BV]).
Suppose $b: \Gamma \rightarrow \ell^2\Gamma$ is a 1-cocycle, as $\Gamma_2$ is nonamenable 
$\ell^2\Gamma$ does not weakly contain the trivial representation for $\Gamma_2$ (see [MV]),
hence $\exists K > 0$, $\gamma_1, \ldots, \gamma_n \in \Gamma_2$ such that $\forall \xi \in \ell^2\Gamma$,
$\| \xi \| \leq K \Sigma_{i = 1}^n \| \lambda(\gamma_i) \xi - \xi \|$.  In particular
we have that $\forall \gamma \in \Gamma_1$, 
$\| b(\gamma) \| \leq K \Sigma_{i = 1}^n \| \lambda(\gamma_i) b(\gamma) - b(\gamma) \|
			= K \Sigma_{i = 1}^n \| \lambda(\gamma) b(\gamma_i) - b(\gamma_i)\|
			\leq 2K \Sigma_{i = 1}^n \| b(\gamma_i) \|$.
Thus we have shown that $b_{| \Gamma_1}$ is bounded and hence we may subtract from $b$ an inner cocycle
and assume that $b_{| \Gamma_1} = 0$.  Therefore we have that $\forall \gamma \in \Gamma_2$, 
$b(\gamma)$ is a $\Gamma_1$-invariant vector, and since $\Gamma_1$ is infinite we must then have that $b(\gamma) = 0$.
Thus we have shown that $b = 0$.

In Theorems 3.3 and 3.5 we will use the same idea as above to show that if $N = Q \overline \otimes B$ is 
a II$_1$ factor where $Q$ is nonamenable and $B$ is diffuse then $N$ must be $L^2$-rigid.  Note that
given a closable derivation $\delta$ on $N$ there is no reason to expect that $Q$ or $B$ is contained 
in the domain of $\delta$, thus it is necessary to use $\tilde \delta_\alpha$ which is everywhere defined and 
by Lemma 2.3 is almost a derivation.  Also note that given $x \in N$, $\tilde \delta_\alpha(x)$ may not be
left and right bounded and thus we have that vectors of the form $y \tilde \delta_\alpha(x)$ and 
$\zeta_\alpha (y) \tilde \delta_\alpha(x)$ may not be close.  This type of situation is handled by Lemma 2.4 and 
using Connes' characterization of amenability [C1] that a factor is amenable if 
the trace has a purification on the minimal tensor product.

Given a free ultrafilter $\omega$, and a unital, tracial, c.p. map $\phi$ on a finite von Neumann algebra $(M, \tau)$
we may extend $\phi$ to a unital, tracial, c.p. map on $M^\omega$ by setting
$\phi(x) = ( \phi(x_n) )_n$ if $x = (x_n)_n$.  If $\{ \phi_\iota \}_\iota$ is a deformation on $M$
which does not converge uniformly on $(M)_1$ 
then the extension to $N^\omega$ does not converge pointwise in $\| \cdot \|_2$ to id.  
We will show however in the next theorem
that if $Q$ is a nonamenable subfactor then not only does an $L^2$-deformation converge pointwise but
it actually converges uniformly to id on $(Q' \cap M^\omega)_1$.

\vskip .05in
\proclaim{3.3. Theorem} Suppose $(N, \tau)$ is a finite von Neumann algebra with normal faithful 
trace $\tau$ and $\{ \eta_\alpha \}_\alpha$ is an $L^2$-deformation for $N$.
If $Q \subset N$ is a nonamenable subfactor and $\omega$ is a free ultrafilter
then $\{ \eta_\alpha \}_\alpha$ converges
uniformly on $(Q' \cap N^\omega)_1$ as $\alpha \rightarrow \infty$.
In particular if $Q \subset N$ is a nonamenable subfactor
then the inclusion $(Q' \cap N \subset N)$ is $L^2$-rigid.
\endproclaim

\vskip .05in
\noindent
{\it Proof}.  Suppose that the deformation 
$\eta_{\alpha}$ does not converge uniformly on $(B)_1$ where $B = Q' \cap N^\omega$.
Then $\exists c > 0$ such that $\forall \alpha > 0$, $\exists x_{\alpha} \in (B)_1$
such that $\| x_{\alpha} - \eta_{\alpha}(x_{\alpha}) \|_2 > c$.
We will show that this implies that $Q$ is amenable.

By [C1] to show that $Q$ is amenable it is enough to show that 
$| \tau(\Sigma a_ib_i^*) | \leq  \| \Sigma a_i \otimes_{\text {min}} b_i^{\text {op}} \|$, 
$\forall  a_1, \ldots, a_n, b_1, \ldots, b_n \in Q_1$.
Note that as a $Q$-$Q$ bimodule $(L^2(M) \overline \otimes L^2(M))^{\oplus \infty}$ is just
a direct sum of coarse correspondences and so the representations
of $Q$ and $Q^{\text {op}}$on $\Cal H$ given by the left and right module structures induce the minimal tensor norm.

Let $\varepsilon > 0$, $a_1, \ldots, a_n, b_1, \ldots, b_n \in Q_1$, since $Q$ is a factor there exists
a finite set $F \subset \Cal U(Q)$ and $0 < \delta \leq \varepsilon$ such that if $\psi \in Q_*$ is a normal state and 
$\| [\psi, u] \| \leq \delta$, $\forall u \in F$ then 
$| \tau( \Sigma a_i b_i^*) - \psi ( \Sigma a_i b_i^*) | < \varepsilon/3$.
Let $F' = F \cup \{ b_i \}_i$, and by Lemma 2.3
let $\alpha_0 > 0$ such that 
$\forall \alpha \geq \alpha_0$, $y \in F'$, $x \in (N)_1$ we have
$\| [\zeta_\alpha (y), \tilde \delta_\alpha(x) ] \| < c^2 \delta/4n + \| \tilde \delta_\alpha([y,x]) \|
	\leq c^2 \delta/4n + 2\| [y, x] \|_2$.

Let $x_\alpha = (x_\alpha^k)_k$ where $\| x_\alpha^k \| \leq 1$, $\forall k \in \Bbb N$
then $\exists k = k(\alpha) \in \Bbb N$ such that 
$\| [y, x_\alpha^k] \|_2 < c^2 \delta/8n$, $\forall y \in F'$,
and $\| x_\alpha^k - \eta_\alpha(x_\alpha^k) \|_2 \geq c$.
Thus if we let $\psi_{\alpha}$ be the state given by 
$z \mapsto \| \tilde \delta_\alpha(x_\alpha^k) \|^{-2} 
	\langle z \tilde \delta_\alpha(x_\alpha^k), \tilde \delta_\alpha(x_\alpha^k) \rangle$,
then by Lemma 2.4 $\exists \alpha \geq \alpha_0$ such that we have 
$| \psi_\alpha \circ \zeta_\alpha (\Sigma a_ib_i^*) 
	- \psi_\alpha( \Sigma \zeta_\alpha(a_i) \zeta_\alpha(b_i^*) ) | < \varepsilon/3$.
Also we may assume that $\forall z \in Q_1$, $u \in F$
$|\psi_{\alpha} \circ \zeta_\alpha (uzu^*) 
	- \psi_\alpha (\zeta_\alpha (uz) \zeta_\alpha(u^*)) | \leq \delta/3$,
and $| \psi_\alpha (\zeta_\alpha (u^*) \zeta_\alpha(uz)) 
	- \psi_{\alpha} \circ \zeta_\alpha (z) | \leq \delta/3$.
Since $F \subset F'$ we have that 
$| \psi_\alpha (\zeta_\alpha (uz) \zeta_\alpha(u^*))
	- \psi_\alpha (\zeta_\alpha (u^*) \zeta_\alpha(uz)) | < \delta/3$,
and so by the triangle inequality we obtain that
$\| [\psi_\alpha \circ \zeta_\alpha, u] \| \leq \delta$, $\forall u \in F$.

Hence 
$$
| \tau ( \Sigma a_i b_i^* ) | \leq | \psi_{\alpha} \circ \zeta_\alpha ( \Sigma a_i b_i^* ) | + \varepsilon/3
$$
$$
\leq | \psi_{\alpha} ( \Sigma \zeta_\alpha(a_i) \zeta_\alpha(b_i^*) ) | 
			+ 2\varepsilon/3
$$
$$
= \| \tilde \delta_\alpha(x_\alpha) \|^{-2} 
	| \langle \Sigma \zeta_\alpha(a_i) \zeta_\alpha(b_i^*) \tilde \delta_\alpha(x_{\alpha}),
					\tilde \delta_\alpha(x_{\alpha}) \rangle | 
			+ 2\varepsilon/3
$$
$$
\leq \| \tilde \delta_\alpha(x_\alpha) \|^{-2} 
	| \langle \Sigma \zeta_\alpha(a_i) \tilde \delta_\alpha(x_{\alpha}) \zeta_\alpha(b_i^*),
					\tilde \delta_\alpha(x_{\alpha}) \rangle | 
			+ \varepsilon
$$
$$
\leq \| \Sigma a_i \otimes_{\text {min}} b_i^{\text {op}} \| + \varepsilon.
$$
Since $\varepsilon$ was arbitrary we have that 
$| \tau ( \Sigma a_i b_i^* ) | \leq \| \Sigma a_i \otimes_{\text {min}} b_i^{\text {op}} \|$
and thus $Q$ is amenable.
\hfill $\square$

\vskip .05in
\proclaim{3.4. Corollary}  Let $\Gamma$ be a countable group and suppose there exists a proper cocycle
$b : \Gamma \rightarrow \ell^2(\Gamma)^{\oplus \infty}$, 
if $B \subset L(\Gamma)$ is diffuse then every 
subfactor of $B' \cap L(\Gamma)$ is amenable.  In particular
all nonamenable subfactors of $L(\Gamma)$ are prime.
\endproclaim

\vskip .05in
\noindent
{\it Proof}.  This follows directly from Theorem 3.3 and remark 3.1.4.
\hfill $\square$

\vskip .05in

We will now show that $L^2$-rigidity passes to normalizers.

\vskip .05in
\proclaim{3.5. Theorem}  Suppose $(N, \tau)$ is a finite von Neumann algebra with normal faithful 
trace $\tau$ and $\{ \eta_\alpha \}_\alpha$ is an $L^2$-deformation for $N$.
If $B \subset N^\omega$ is a
diffuse von Neumann subalgebra 
such that $\{ \eta_\alpha \}_\alpha$ converges uniformly on $(B)_1$,
then $\{ \eta_\alpha \}_\alpha$ converges uniformly 
on $W^*(N \cap \Cal N_{N^\omega}(B))_1$.
In particular if $B \subset N$ is a diffuse von Neumann subalgebra and
$(B \subset N)$ is $L^2$-rigid, then $(W^*(\Cal N_N(B)) \subset N)$ is also $L^2$-rigid.
\endproclaim

\vskip .05in
\noindent
{\it Proof}.  Let $1 \geq \varepsilon > 0$, using Lemma 2.3 
$\exists \alpha_0' > 0$ such that $\forall \alpha > \alpha_0'$, $x = (x_n)_n \in B_1$ (with $\| x_n \| \leq 1$), 
$y \in N_1$ we have 
$\lim_{n \rightarrow \omega} \| \eta_{\alpha}(x_n) - x_n \|_2 < \varepsilon / 8$, and
$\lim_{n \rightarrow \omega} \| \zeta_\alpha (x_n) \tilde \delta_\alpha (y) + \tilde \delta_\alpha(x_n)\zeta_\alpha(y) 
			 - \tilde \delta_\alpha(x_ny) \| < \varepsilon/8$.
Take $v \in N \cap \Cal N_{N^\omega}(B)$ and $\alpha > \alpha_0'$, then  
since $B$ is diffuse, by the mixing property
of the coarse correspondence we have that
$\exists u = (u_n)_n \in \Cal U(B)$ (with $u_n \in \Cal U(N)$) such that 
$\| \tilde \delta_\alpha( v ) \|_2
  \leq \lim_{n \rightarrow \omega} 2\| \zeta_\alpha(u_n) \tilde \delta_\alpha( v ) \zeta_\alpha (v^*u_n^*v) 
			- \tilde \delta_\alpha( v ) \|_2$.
Hence we have:
$$
\| v - \eta_\alpha(v) \|_2^2 \leq \| \tilde \delta_\alpha( v ) \|^2
$$
$$
\leq  \lim_{n \rightarrow \omega} 4\| \zeta_\alpha(u_n) \tilde \delta_\alpha( v ) \zeta_\alpha(v^*u_n^*v)
			- \tilde \delta_\alpha( v )  \|^2
$$
$$
\leq \lim_{n \rightarrow \omega} 4( \| \tilde \delta_\alpha(u_n) \| 
	+ \| \tilde \delta_\alpha(u_n) \zeta_\alpha(v) 
		+ \zeta_\alpha(u_n) \tilde \delta_\alpha(v) - \tilde \delta_\alpha(u_nv) \|
$$
$$
	+ \| \tilde \delta_\alpha(v^*u_n^*v) \|
	+ \| \tilde \delta_\alpha(u_nv) \zeta_\alpha (v^*u_n^*v) 
		+ \zeta_\alpha(u_nv) \tilde \delta_\alpha(v^*u_n^*v) - \tilde \delta_\alpha(v) \| )^2
	< \varepsilon^2,
$$ 
as the maps $\eta_\alpha$ are tracial the result then follows by standard arguments (see [P2]).
\hfill $\square$

\vskip .05in
\proclaim{3.6. Corollary}  If $N$ is a nonamenable II$_1$ factor which is non-prime
or has property $\Gamma$, then $N$ is $L^2$-rigid.
\endproclaim

\vskip .05in
\noindent
{\it Proof}.  If $N = Q \overline \otimes B$ with $Q$ a nonamenable factor then by
Theorem 3.3 we have that $(B \subset N)$ is $L^2$-rigid.  If $B$ is diffuse then
by Theorem 3.5 we then have that $N$ is $L^2$-rigid.

Also if $N$ is a nonamenable factor then by Theorem 3.3
if $\omega$ is a free ultrafilter then any $L^2$-deformation converges uniformly on 
$(N' \cap N^\omega)_1$, if $N$ has property $\Gamma$ then $N' \cap N^\omega$ is 
diffuse and so from Theorem 3.5 we would have that the $L^2$-deformation
converges uniformly on $(N)_1$.
\hfill $\square$

\vskip .05in
\proclaim{3.7. Corollary}
Let $N$ be a finite von Neumann algebra such that $N$ is a free product
of diffuse finite von Neumann algebras or 
let $N = L\Gamma$ where $\Gamma$ is a countable group with $H^1(\Gamma, \ell^2(\Gamma)) \not= \{ 0 \}$.

\noindent
$1.$ If $B \subset N$ is a regular diffuse subalgebra then every subfactor of $B' \cap N$ is amenable.

\noindent
$2.$ Any nonamenable regular subfactor of $N$ is prime and does not have properties $\Gamma$ or (T).
\endproclaim

\vskip .05in
\noindent
{\it Proof}.  $1.$ If $Q \subset B' \cap N$ is a nonamenable subfactor 
then by Theorem 3.3 we would have that $(B \subset N)$ 
is $L^2$-rigid, hence by Theorem 3.5 we would have that $N$ is $L^2$-rigid and thus the result follows from
remark 3.1.3.

\vskip .05in
\noindent
$2.$  By Corollary 3.6 and Theorem 3.5 if $N$ has a regular subfactor 
which is non-prime or has properties $\Gamma$ or (T)
then $N$ is $L^2$-rigid and so as above the result follows from remark 3.1.3.
\hfill $\square$

\vskip .05in
Note that if $\Gamma$ is finitely generated and non-amenable then by [BV]
$H^1(\Gamma, \ell^2(\Gamma)) \not= \{ 0 \}$
if and only if $b_1^{(2)}(\Gamma) > 0$.  For nonamenable groups which are not finitely generated it follows from
a result of Gaboriau that if $H^1(\Gamma, \ell^2(\Gamma)) \not= \{ 0 \}$ then $b_1^{(2)}(\Gamma) > 0$,
however the reverse implication is open [MV].

\heading 4. $L^2$-rigid subalgebras in free product factors. \endheading

Let $(M_i, \tau_i)$, $i = 1,2$ be finite von Neumann algebras, denote $M = M_1 * M_2$.
Let $\delta_i:M_1 *_{\text {Alg}} M_2 \rightarrow L^2(M) \otimes L^2(M)$ be the unique derivation
which satisfies $\delta_i(x) = x \otimes 1 - 1 \otimes x$, $\forall x \in M_i$ and 
$\delta_i(y) = 0$, $\forall y \in M_j$ where $j \not=i$.  Then as above we have that
$\phi_s^1 = (e^{-2s}{\text {id}} + (1-e^{-2s})\tau) * {\text {id}}$,
and $\phi_s^2 = {\text {id}} * (e^{-2s}{\text {id}} + (1-e^{-2s})\tau)$ are the associated semigroups
of c.p. maps. 

If $Q$ is an $L^2$-rigid subalgebra of $M$ then we may interpret the fact that 
the above deformations converge uniformly on $(Q)_1$
as saying that $Q$ has ``bounded word length''.  Thus one would expect that a corner of $Q$ embeds
into either $M_1$ or $M_2$.  We will show in this section that this is indeed the case, we do this by first 
showing that $Q$ must be rigid with respect to the deformations used in [IPP], then we may apply the
word reduction argument in [IPP] (Theorem 4.3) which gives the result.

Recall that if we let $\Cal H_i^0 = L^2(M_i) \ominus \Bbb C$ then we may decompose $L^2(M_1 * M_2)$ in the 
usual way as
$$
L^2(M_1 * M_2) = \Bbb C \oplus \bigoplus_{n \geq 1} 
		\bigoplus_{{ i_j \in \{ 1, 2 \} } \atop {i_1 \not= i_2, \cdots, i_{n-1} \not= i_n }}
			\Cal H_{i_1}^0 \otimes \Cal H_{i_2}^0 \otimes \cdots \otimes \Cal H_{i_n}^0.
$$

\vskip .05in
\proclaim{4.1. Lemma}  Let $(M_1, \tau_1)$, $(M_2, \tau_2)$ be finite von Neumann algebras.  
As in 2.2.2 of [IPP] denote 
$M = M_1 * M_2$, $\tilde M_j = M_j * L(\Bbb Z)$, $j = 1,2$, and 
$\tilde M = \tilde M_1 * \tilde M_2 = M * L(\Bbb F_2)$.  
Let $h_j \in L(\Bbb F_2)$ be self-adjoint elements such that
$u_j = \exp(\pi i h_j)$, where $u_1, u_2 \in L(\Bbb F_2)$
are the canonical generators of $L(\Bbb F_2)$.  Let $u_j^t = \exp(\pi i t h_j)$, and set
$\theta_t = {\text {Ad}}(u_1^t) * {\text {Ad}}(u_2^t) \in {\text {Aut}}(\tilde M)$, a one parameter group of 
automorphisms.  Suppose $Q \subset M$ is a von Neumann subalgebra, then
the deformation $\{ \theta_t \}_t$ converges uniformly on $(Q)_1$ as $t \rightarrow 0$ 
if and only if the deformations
$\{ \phi_s^j \}_s$ converge uniformly on $(Q)_1$ as $s \rightarrow 0$, $j = 1,2$.
\endproclaim

\vskip .05in
\noindent
{\it Proof}.  Let $\varepsilon_0 > 0$ such that $\tau(u_j^t) \not= 0$, $\forall t < \varepsilon_0$, $j = 1,2$.
Let $t < \varepsilon_0$, it is then a simple exercise to check that if $f_j(t) = - \log ( | \tau(u_j^t) | )$
then $\tau ( \theta_t (x)x^* ) = \tau ( \phi_{f_j(t)}^j (x) x^*)$, $\forall x \in M_j$.  
In fact using the direct sum decomposition above one sees that if 
$x = x_1x_2 \cdots x_n$, where $i_j \in \{ 1, 2 \}$, $j \leq n$, $i_1 \not= i_2, \cdots, i_{n-1} \not= i_n$,
and $x_j \in \Cal H_{i_j}^0$, $\forall j \leq n$.
Then in fact we have that
$\tau ( \theta_t (x)x^* ) = \tau ( \theta_t (x_1) x_1^* ) \cdots \tau ( \theta_t (x_n) x_n^* )
	= \tau ( \phi_{f_{i_1}(t)}^{i_1} (x_1) x_1^* ) \cdots  \tau ( \phi_{f_{i_n}(t)}^{i_n} (x_n) x_n^* ) 
	= \tau ( \phi_{f_1(t)}^1 \circ \phi_{f_2(t)}^2 (x) x^* )$.

Moreover since both of the maps ${\theta_t}_{|M}$ and $\phi_{f_1(t)}^1 \circ \phi_{f_2(t)}^2$ take orthogonal vectors 
to orthogonal vectors we have that
$\tau ( \theta_t (x) x^*) = \tau ( \phi_{f_1(t)}^1 \circ \phi_{f_2(t)}^2 (x) x^* )$, $\forall x \in M$.

Since $\| \phi_{f_1(t)}^1 \circ \phi_{f_2(t)}^2 (x) - x \|_2 \geq \| \phi_{f_j(t)}^j (x) - x \|_2$,
$\forall x \in M$, $j = 1,2$, and since
$f_j(t) \rightarrow 0$, $j = 1,2$ as $t \rightarrow 0$ the result follows easily.
\hfill $\square$

\vskip .05in
\proclaim{4.2. Corollary}  Let $M_1$ and $M_2$ be separable $II_1$ factors, and let $M = M_1 * M_2$.  
If $(Q \subset M)$ is $L^2$-rigid then there exists a unique pair of projections
$q_1, q_2 \in Q' \cap M$ such that $q_1 + q_2 = 1$, and 
$u_i(Qq_i)u_i^* \subset M_i$ for some unitaries $u_i \in \Cal U(M)$, $i = 1,2$.
Moreover, these projections lie in the center of $Q' \cap M$.
\endproclaim

\vskip .05in
\noindent
{\it Proof}.  
Suppose $(Q \subset M)$ is $L^2$-rigid, then by definition we have that the deformations $\{ \phi_s^j \}_s$, 
converge uniformly to id on $(Q)_1$ as $s \rightarrow 0$, hence by Lemma 4.1 the deformation
$\{ \theta_t \}_t$ also converges uniformly on $(Q)_1$ as $t \rightarrow 0$.  A check of 
Theorem 4.3 in [IPP] shows that these are the only two facts used from the rigid inclusion.
Thus the result follows from Theorems 4.3 and 5.1 in [IPP].
\hfill $\square$

\heading 5. Unique prime factorization and non-$L^2$-rigid factors. \endheading

In this section we will adapt Theorems 3.3 and 3.5 and use Popa'a intertwining technique along with
the results in [OP] in order to show 
that if $M_i$ are II$_1$ factors which have derivations into $L^2(M_i) \otimes L^2(M_i)$
which do not ``vanish'' then the tensor product has unique prime factorization (up to amplification
and unitary conjugation of the factors).  In order to satisfy the conditions of 
Popa's intertwining criteria (Theorem 2.1 in [P5])
it will be necessary to assume that the derivation is actually densely defined on $M_i$. 
This is a formally stronger condition then the negation of $L^2$-rigidity, however note that both 
examples 1 and 2 in section 1 satisfy this condition.  
For the following theorem if $M = M_1 \overline \otimes M_2 \overline \otimes \cdots \overline \otimes M_m$
then we will denote by $M_i'$ the resulting von Neumann subalgebra obtained by replacing
$M_i$ with $\Bbb C1$ so that $M = M_i \overline \otimes M_i'$,

\proclaim{5.1. Theorem}  Let $M_i$ be nonamenable II$_1$ factors $1 \leq i \leq m$,
suppose that each $M_i$ has a densely defined real closable derivation
into $(L^2(M_i) \otimes L^2(M_i))^{\oplus \infty}$ such that the 
associated $L^2$-deformation does not converge uniformly on $(M_i)_1$. 
Let $M = M_1 \overline \otimes M_2 \overline \otimes \cdots \overline \otimes M_m$.
Assume that $B \subset M$ is a regular type $II_1$ factor such that $B' \cap M$ is 
a nonamenable subfactor.  Then $\exists k \in \{ 1, \ldots, m \}$, $t > 0$
and a unitary element $u \in \Cal U(M)$ such that 
$uBu^* \subset (M_k')^t \otimes \Bbb C \subset (M_k')^t \overline \otimes (M_k)^{1/t} = M$.
If in addition we have that the $L^2$-deformations above may all be taken compact then $B$ need not
be regular.
\endproclaim

\vskip .05in
\noindent
{\it Proof}.  Let $\delta_0^i: M_i \rightarrow (L^2(M_i) \overline \otimes L^2(M_i))^{\oplus \infty}$ be a 
densely defined closable real derivation such that the corresponding deformation $\{ \eta_\alpha^i \}$ 
does not converge uniformly on $(M_i)_1$.  
Then we may embed $(L^2(M_i) \overline \otimes L^2(M_i))^{\oplus \infty}$
into $\Cal H_i = (L^2(M) \overline \otimes_{M_i'} L^2(M) )^{\oplus \infty}$ in the natural way 
as $M_i$-$M_i$ Hilbert bimodules and we then may extend $\delta_0^i$ to 
a densely defined closable real derivation $\delta^i$ on $M$ by setting 
$\delta^i(x) = 0$, $\forall x \in M_i'$.  We denote by $\{ \hat \eta_\alpha^i \}$
the corresponding deformations on $M$, so that $\hat \eta_\alpha^i = \eta_\alpha^i \otimes$ id,
also let $\hat \zeta_\alpha^i = \zeta_\alpha^i \otimes$ id $= (\hat \eta_\alpha^i)^{1/2}$.

We will proceed as in Theorem 3.3 to show that if each $\{ \hat \eta_\alpha^i \}$ 
does not converge uniformly on $(B)_1$ then 
we must have that $Q = B' \cap M$ is amenable.  Indeed if this is the case then $\exists c > 0$, such that
$\forall \alpha > 0$, $i \leq m$, $\exists x_\alpha^i \in (B)_1$ such 
that $\| x_\alpha^i - \hat \eta_\alpha^i(x_\alpha^i) \|_2 \geq c$.
Let $\varepsilon > 0$, $a_1, \ldots, a_n$, $b_1, \ldots, b_n \in Q_1$.
Let $F \subset \Cal U(Q)$ finite, and $0 < \delta < \varepsilon$ such that 
if $\psi \in Q_*$ is a normal state and 
$\| [ \psi, u ] \| \leq \delta$, $\forall u \in F$ then 
$| \tau (\Sigma a_i b_i^*) - \psi( \Sigma a_i b_i^*) | < \varepsilon / 3$.
Let $F' = F \cup \{ b_i \}_i$, by Lemma 2.3 let $\alpha_0 > 0$ such that $\forall \alpha \geq \alpha_0$,
$y \in F'$, $x \in (B)_1$ we have $\| [\hat \zeta_\alpha(y), \tilde \delta_\alpha^i(x) ] \| < c^2\delta/4nm$
and $\| [ \check \zeta_\alpha^i (y), x] \|_2 < c^2\delta/8nm$
where
$\hat \zeta_\alpha = \zeta_\alpha^1 \circ \cdots \circ \zeta_\alpha^m$
and $\check \zeta_\alpha^i$ is obtained by omitting $\zeta_\alpha^i$ from $\hat \zeta_\alpha$.  

Let $_M\Cal H_M = \Cal H_1 \overline \otimes_M \Cal H_2 \overline \otimes_M \cdots \overline \otimes_M \Cal H_m$,
and note that $\Cal H$ may be embedded into $(L^2(M) \overline \otimes L^2(M) )^{\oplus \infty}$
as $M$-$M$ Hilbert bimodules.
Let $\xi_\alpha = \tilde \delta_\alpha^1(x_\alpha^1) \otimes \cdots \otimes \tilde \delta_\alpha^m(x_\alpha^m) \in \Cal H$,
and let $\psi_\alpha$
be the normal state given by 
$z \mapsto \| \xi_\alpha \|^{-2} \langle z \xi_\alpha, \xi_\alpha \rangle$.
Then $\psi_\alpha = \psi_\alpha^1 \otimes \cdots \otimes \psi_\alpha^m$ where $\psi_\alpha^i$ is the 
state given by 
$z \mapsto \| \tilde \delta_\alpha^i(x_\alpha^i) \|^{-2} 
	\langle z \tilde \delta_\alpha^i(x_\alpha^i), \tilde \delta_\alpha^i(x_\alpha^i) \rangle$
hence we may apply Lemma 2.4 (parts (ii) and (iii)) to insure that for large enough $\alpha$ we have
$| \tau( \Sigma a_ib_i^* ) - \psi_\alpha \circ \hat \zeta_\alpha ( \Sigma a_ib_i^* ) | < \varepsilon/3$,
and $| \psi_\alpha \circ \hat \zeta_\alpha ( \Sigma a_ib_i^* )
	- \psi_\alpha( \Sigma \hat \zeta_\alpha (a_i) \hat \zeta_\alpha (b_i^*) ) | < \varepsilon/3$.
Then the same proof in 3.3 shows that we obtain
$| \tau(\Sigma a_i b_i^* ) | \leq \| \Sigma a_i \otimes_{\text {min}} b_i^{\text {op}} \| + \varepsilon$.
Since $\varepsilon$ was arbitrary we have that $Q$ is amenable.

Therefore if $B' \cap M$ is a nonamenable factor then we have shown that
$\exists k \leq m$ such that the deformation $\{ \hat \eta_\alpha^k \}$ converges uniformly 
on $(B)_1$.
Next we show that if this is the case then we have that a corner of $B$ embeds into
$M_k'$ inside of $M$, i.e. there exists a non-zero projection $f$ in 
$B' \cap \langle M, e_{M_k'} \rangle$ of finite trace $Tr = Tr_{\langle M, E_{M_k'} \rangle}$.

If we do not have that a corner of $B$ embeds into $M_k'$ inside of $M$ then by Corollary 2.3 of [P5]
there exists a sequence of unitaries $\{ u_n \}_n \subset \Cal U(B)$ such that
$\forall x \in M$, $\| E_{M_k'}(x u_n ) \|_2 \rightarrow 0$, as $n \rightarrow \infty$.
Since $\hat {\zeta_\alpha^k}_{|M_k'} =$ id we have that 
$\forall x \in M$, $\| E_{M_k'}(x \hat \zeta_\alpha^k(u_n) ) \|_2 \rightarrow 0$, as $n \rightarrow \infty$,
and since $M_k'$ is regular in $M$ this implies 
$\| E_{M_k'}( x \hat \zeta_\alpha^k(u_n) y ) \|_2 \rightarrow 0$, as $n \rightarrow \infty$, $\forall x,y \in M$.
In particular this shows that $\forall v \in \Cal N_M(B)$, $\exists u \in \Cal U(B)$ such that
$$
\| \zeta_\alpha^k(u) \tilde \delta_\alpha^k(v) \zeta_\alpha^k(v^*u^*v) 
		- \tilde \delta_\alpha^k(v) \| 
	\geq \| \tilde \delta_\alpha^k(v) \| \tag 5.1.1
$$

On the other hand since $B$ is regular and since 
$\{ \eta_\alpha^k \}_\alpha$ does not converge uniformly on $(M)_1$,
$\exists c_0 > 0$ such that $\forall \alpha > 0$, $\exists v_\alpha \in \Cal N_M(B)$ such that
$\| \tilde \delta_\alpha^k(v_\alpha) \| \geq \|v_\alpha - \eta_\alpha^k(v_\alpha) \|_2 \geq c_0$.
By Lemma 2.3 $\forall \varepsilon > 0$, $\exists \alpha_0 > 0$ such that
$\forall \alpha \geq \alpha_0$, $u \in \Cal U(B)$ we have that
$$
\| \zeta_\alpha^k(u) \tilde \delta_\alpha^k(v_\alpha) \zeta_\alpha^k(v_\alpha^*u^*v_\alpha) 
	- \tilde \delta_\alpha^k(v_\alpha) \|
< \varepsilon.
$$

Thus for $\varepsilon < c_0$ we have
$$
\| \zeta_\alpha^k(u) \tilde \delta_\alpha^k(v_\alpha) \zeta_\alpha^k(v_\alpha^*u^*v_\alpha) 
	- \tilde \delta_\alpha^k(v_\alpha) \|
< \| \tilde \delta_\alpha^k(v_\alpha) \|,
$$
for each $u \in \Cal U(B)$, which contradicts $(5.1.1)$.

If $B$ is not regular but each deformation is compact then we may apply the proof of Theorem 6.2 in [P4]
to show that a corner of $B$ embeds into $M_k'$ inside of $M$ in this case also.

Thus in either case since
$B' \cap M$ is a factor we may then apply Proposition 12 in [OP] to obtain the result.
\hfill $\square$

\vskip .05in
As a consequence of the previous theorem, 
we obtain from [OP] the following unique prime factorization result.

\vskip .05in
\proclaim{5.2. Corollary}  Let $M_i$ be nonamenable II$_1$ factors $1 \leq i \leq m$,
suppose that each $M_i$ has a densely defined real closable derivation
into $(L^2(M_i) \otimes L^2(M_i))^{\oplus \infty}$ such that the 
associated $L^2$-deformation does not converge uniformly on $(M_i)_1$.  
Assume $N_1 \overline \otimes \cdots \overline \otimes N_n = M_1 \overline \otimes \cdots \overline \otimes M_m$,
for some prime II$_1$ factors $N_1, \ldots, N_n$, then $n = m$ and there exist
$t_1, t_2, \ldots, t_m > 0$ with $t_1t_2 \cdots t_n = 1$ such that after permutation of indices 
and unitary conjugacy we have $N_k^{t_k} = M_k$, $\forall k \leq m$.
\endproclaim

\heading  References\endheading

\item{[BV]} M.E.B. Bekka, A. Valette: {\it Group cohomology, harmonic functions and the first $L^2$-Betti number},
Potential Analysis, {\bf 6} (1997), 313-326.

\item{[CiS]} F. Cipriani, J.-L. Sauvageot: {\it Derivations as square roots of Dirichlet forms},
J. of Functional Analysis {\bf 201} (2003), 78-120.

\item{[C1]} A. Connes: {\it Classification of injective factors}, 
Ann. of Math., {\bf 104} (1976), 73-115.

\item{[C2]} A. Connes: {\it A type II$_1$ factor with countable fundamental group},
J. Operator Theory, {\bf 4} (1980), 151-153.

\item{[C3]} A. Connes: {\it Classification des facteurs}, Proc. Symp. Pure Math., {\bf 38}
(Amer. Math. Soc., 1982), 43-109.

\item{[CJ]} A. Connes, V.F.R. Jones: {\it Property (T) for von Neumann algebras},
Bull. London Math. Soc., {\bf 17} (1985), 57-62.

\item{[CSh]} A. Connes, D. Shlyakhtenko: {\it $L^2$-homology for von Neumann algebras}, 
J. Reine Angew. Math. 586 (2005), 125-168.

\item{[DL]} E.B. Davies, J.M. Lindsay: {\it Non-commutative symmetric Markov semigroups},
Math. Zeitschrift {\bf 210} (1992), 379-411.

\item{[Ge]} L. Ge: {\it Applications of free entropy to finite von Neumann algebras. II}, 
Ann. of Math. (2) {\bf 147} (1998), no. 1, 143-157. 

\item{[H]} U. Haagerup: {\it An example of a nonnuclear $C^*$-algebra, which has the metric approximation property},
Invent. Math. {\bf 50} (1979), 279-293.

\item{[IPP]} A. Ioana, J. Peterson, S. Popa: {\it Amalgamated free products of $w$-rigid factors and calculation of their symmetry groups},
Preprint 2005, math.OA/0505589.

\item{[J]} K. Jung: {\it Strongly } $1$-{\it bounded von Neumann algebras},
Preprint 2005. math.OA/0510576, to appear in GAFA.

\item{[MR]} Z.-M. Ma, M. R\" ockner: ``Introduction to the theory of (non-symmetric) Dirichlet forms'',
Universitext. Springer, Berlin. 1992.

\item{[MV]} F. Martin, A. Valette: {\it On the first $L^p$-cohomology of discrete groups}, preprint 2006.

\item{[MvN]} F. J. Murray, J. von Neumann: {\it On rings of operators IV}, 
Ann. of Math. (2) {\bf 44} (1943), 716-808.

\item{[O1]} N. Ozawa: {\it Solid von Neumann algebras}, Acta Math. {\bf 192} (2004), no. 1, 111-117.

\item{[O2]} N. Ozawa: {\it A Kurosh type theorem for type II$\sb 1$ factors}, Preprint 2004, math.OA/0401121.

\item{[OP]} N. Ozawa, S. Popa: {\it Some prime factorization results for II$\sb 1$ factors},
Invent. Math. {\bf 156} (2004), 223-234.

\item{[Pe]} J. Peterson: {\it A} $1$-{\it cohomology characterization of property (T) in von Neumann algebras}, Preprint 2004, 
math.OA/0409527

\item{[P1]} S. Popa: {\it Orthogonal pairs of $*$-subalgebras in finite von Neumann algebras}, 
J. Operator Theory {\bf 9}, no. 2, 253-268.

\item{[P2]} S. Popa: {\it Correspondences}, INCREST preprint 1986,
unpublished.

\item{[P3]} S. Popa: {\it Some rigidity results for non-commutative Bernoulli shifts},
J. Funct. Anal., {\bf 230} (2006), no.2, 273-328.

\item{[P4]} S. Popa: {\it On a class of type II$_1$ factors with Betti numbers invariants},
MSRI Preprint 2001-005, math.OA/0209130, to appear in Annals of Mathematics.

\item{[P5]} S. Popa: {\it Strong rigidity of II$\sb 1$ factors arising from malleable actions of $w$-rigid groups I, II},
preprints 2003 and 2004, math.OA/0305306 and math.OA/0407103.

\item{[S1]} J.-L. Sauvageot: {\it Tangent bimodules and locality for dissipative operators on $C^*$-algebras}, 
Quantum probability and appl.,
IV, Lecture notes in Math. {\bf 1396} (1989), 322-338.

\item{[S2]} J.-L. Sauvageot: {\it Quantum Dirichlet forms, differential calculus and semigroups}, 
Quantum probability and appl., V, Lecture
notes in Math. {\bf 1442} (1990), 334-346.

\item{[S3]} J.-L. Sauvageot: {\it Strong Feller semigroups on $C^*$-algebras}, 
J. Operator Theory {\bf 42} (1999), 83-102.

\item{[T]} A. Thom: {\it $L^2$-cohomology for von Neumann algebras}, Preprint 2006,
math.OA/0601447.

\item{[V]} D. Voiculescu: {\it The analogues of entropy and of Fisher's 
information measure in free probability theory, V},  Invent. Math. {\bf 132} (1998), 189-227.
\enddocument